\documentclass[12pt, reqno]{amsart}
\usepackage{amsmath, amsthm, amscd, amsfonts, amssymb, graphicx, hyperref,color}
\newtheorem{theorem}{Theorem}[section]
\newtheorem{lemma}[theorem]{Lemma}

\newtheorem{corollary}[theorem]{Corollary}
\theoremstyle{definition}

\theoremstyle{remark}
\newtheorem{remark}[theorem]{Remark}
\numberwithin{equation}{section}
\begin{document}
%%%%%%%%%%%%%%%%%%%%%%%%%%%%%%%%%%%%%%%%%%%%%%%%%%
%%%%%%%%%%%%%%%%%%%%%%%%%%%%%%%%%%%%%%%%%%%%%%%%%
\setcounter{page}{1}
%%%%%%%%%%%%%%%%%%%%%%%%%%%%%%%%%%%%%%%%%%%%%

%-------------------------- Pleased do not change the following line-------------------------------------------
%\noindent \textcolor[rgb]{0.99,0.00,0.00}{}\\[.5in]
%--------------------------------------------------------------------------------------------------------------

\title[A perturbation result  of m-accretive linear operators]
      {A perturbation result  of m-accretive linear operators in Hilbert spaces}
       
%%%%%%%%%%%%%%%%%%%%%%%%%%%%%%%%%%%%%%%%%%%%%%%%%
\author[ M. Benharrat]{Mohammed Benharrat$^{1*}$}
%%%%%%%%%%%%%%%%%%%%%%%%%%%%%%%%%%%%%%%%%%%%%%%
%%%%%%%%%%%%%%%%%%%%%%%%%%%%%%%%%%%%%%%%%%%%%%%%
\address{$^{1}$ D\'{e}partement de G\'{e}nie des syst\'{e}mes,
	Ecole Nationale Polytechnique d'Oran-Maurice Audin (Ex. ENSET d'Oran), 
	BP 1523 Oran-El M'naouar, 31000 Oran, Alg\'{e}rie.}
\email{\textcolor[rgb]{0.00,0.00,0.84}{mohammed.benharrat@enp-oran.dz, mohammed.benharrat@gmail.com}}

%%%%%%%%%%%%%%%%%%%%%%%%%%%%%%%%%%%%%%%%%%%%%%%%%%%%%%%%%%%%

%\dedicatory{This paper is dedicated to Professor ABCD}



\date{25/07/2020.
\newline \indent $^{*}$ Corresponding author\\
This work was supported by the Algerian research project: PRFU, no. C00L03ES310120180002.}
%%%%%%%%%%%%%%%%%%%%%%%%%%%%%%%%%%%%%%%%%%%%%%%%%%%%%%%%%%%%%%%%
\begin{abstract}
%%%%%%%%%%%%%%%%%%%%%%%%%%%%%%%%%%%%%%%%%%%%%%%%%%%%%%%%%%%%%
A new  sufficient condition is given for the sum  of  linear m-accretive operator and accretive operator one  in a Hilbert space to be m-accretive. As an application, an extended result to the operator-norm error bound estimate for the exponential Trotter-Kato product formula is given. 
\end{abstract}
%%%%%%%%%%%%%%%%%%%%%%%%%%%%%%%%%%%%%%%%
%%%%%%%%%%%%%%%%%%%%%%%%%%%%%%%%%%%%%%%%%%%%%%%%%%%%%%%%%%%%%%%
 \maketitle
%%%%%%%%%%%%%%%%%%%%%%%%%%%%%%%%%%%%%%%%%%%%%%%%%%%%%%%%%%%%%%
%%%%%%%%%%%%%%%%%%%%%%%%%%%%%%%%%%%%%%%%%%%
%Introduction
%%%%%%%%%%%%%%%%%%%%%%%%%%%%%%%%%%%%%%%%%%

\section{introduction}
%%%%%%%%%%%%%%%%%%%%%%%%%%%%%%%%%%
%Let $T$ be an unbounded linear operator with domain $\mathcal{D} (T)$ in a complex Hilbert space $\mathcal{H}$. The aim of  this paper is to discuss known results and establish new perturbation results on the m-accretivity of the operator $S = T + A$, where $A$ is a non-symmetric operator which is subordinate to the original T in some sense. We
%will consider various types of subordination, and at the end of the paper we will
%look at concrete, relatively simple examples, discuss potential applications for the
%results obtained, and carry out a comparative analysis of different methods.

%Let be $\mathcal{H}$ a complex  separable Hilbert space with inner product $<\cdot, \cdot>$ and norm $\|\cdot\|$.  Given a linear operator $T$ on $\mathcal{H}$,  we denote by  $\mathcal{D} (T)$, $\mathcal{N}(T)$, and $\mathcal{R}(T)$ the domain, the kernel and the range of $T$ respectively. We also denote by   $\sigma (T)=\mathbb{C}\backslash \rho (T)$ and $\rho (T)$ the spectrum and the resolvent set of  $T$, where $\rho (A)=\left\{ \lambda \in \mathbb{C}: \quad  \lambda I-T: \mathcal{D} (T) \longrightarrow \mathcal{H}  \text{ is bijective } \right\} $ and $I$ is the operator identity on $\mathcal{H}$. For $%
%\lambda \in \rho (T),$ the inverse $(\lambda I-T)^{-1}$ is, by the closed graph theorem, a bounded operator on $\mathcal{H}$ and will be called the
%resolvent of $T$ at the point $\lambda$. We now recall the definition of the  accretive operators. 
%%%%%%%%%%%%%%%%%%%%%%%%%%%%%%%%%%%%%%%

A linear operator $T$ with domain $\mathcal{D}(T)$ in a complex Hilbert space $\mathcal{H}$ is said to be accretive if
$$
Re<Tx,x>\geq 0 \qquad \text{ for all } x\in \mathcal{D}(T)
$$
or, equivalently if
$$\Vert(\lambda +T)x\Vert\geq \lambda \Vert x\Vert \qquad \text{ for all } x\in \mathcal{D}(T) \text{ and  }\lambda>0.$$
Further, if $\mathcal{R}(\lambda +T)=\mathcal{H}$ for some (and hence for every) $\lambda>0$, we say
that $T$ is m-accretive. In particular, every m-accretive operator is accretive
and  closed densely defined, its adjoint is also m-accretive (cf. \cite{Kato}, p. 279). Furthermore, $$ (\lambda + T)^{-1} \in \mathcal{B} \ (\mathcal{H}) \quad  \text{ and } \quad \left\|(\lambda + T)^{-1}\right\| \leq \frac{1}{\lambda}   \text{ for } \lambda >0,$$
where, $\mathcal{B} (\mathcal{H})$  denote the Banach space of all bounded linear operators on $\mathcal{H}$. In particular, a bounded accretive  operator  is m-accretive. 

Consider two linear operators $T$ and $A$ in the Hilbert space $\mathcal{H}$, such that $\mathcal{D} (T)\subset \mathcal{D} (A)$. Assume furthermore that $T$ is m-accretive and   $A$ is an accretive operator. Then the question is:

Under which conditions the sum $T + B$ is m-accretive?

% We %%%%%%%%%%%%%%%%%%%%%%%%%%%%%%%%%%%%%%%%%%%%%%%%%%%%%

%define the sum operator by
%$$Sx = Tx + Ax, \quad \text{ for } x\in  \mathcal{D} (S)=\mathcal{D} (T)\subset \mathcal{D} (A).$$
%We exclude from our consideration the trivial case $\mathcal{D} (T)\subset \mathcal{D} (A)=\{0\}.$

Many papers have been devoted to this problem and most results treat pairs $T$, $A$ of relatively bounded or
resolvent commuting operators. We refer the reader to \cite{Chernoff72, Engel95, HessK70, Krol2009, Okazawa77, Okazawa2002, Sobajima14, Wust1971, Yoshikawa72, Yosida65}.  Since $T$ is closed it follows that
there are two nonnegative constants $a$, $b$ such that
\begin{equation}\label{Tbounded}
\left\|Ax\right\|^2\leq a \left\|x\right\|^2+b\left\|Tx\right\|^2  ,\quad  \text{ for all }  x\in \mathcal{D}  (T)\subset  \mathcal{D}  (A).
\end{equation}
In this case, $A$  is  called relatively bounded with respect to $T$ or simply $T$-bounded, and refer to $b$ as a relative bound.
%We say that B is bounded relative to A
%An operator $T$ is called dissipative (resp. m-dissipative) if  $-T$ accretive (resp. m-accretive). 
%A normal operator $T$ (bounded or not) is m-accretive if and only if its spectrum is contained in the
%half-plane $Re(z)\geq 0$ of the plane of the complex numbers $z$. Hence a normal accretive  operator  is m-accretive.
Gustafson \cite{Gustafson66}, generalizing basic work of Rellich, Kato, and others
(cf. \cite{Kato}), showed that that $T+A$ is also m-accretive
if $A$ is $T$-bounded, with $b <1$ (see  \cite[Theorem 2.]{Gustafson66}).  
Okazawa showed in \cite{Okazawa69} that the closure of the sum $T+A$ is m-accretive, if  the bounded  operator $ A(t+ T)^{-1}$ on $\mathcal{H}$ is a contraction for some $t>0$, \cite[Theorem 1.]{Okazawa69}. In particular, he also showed that the validity of \eqref{Tbounded} with $b=1$ implies that the closure of $T+A$ is m-accretive, \cite[Corollary 1.]{Okazawa69}.  Later, the same author in  \cite{Okazawa73} gave a variant of perturbation  by  assumed  the existence of two nonnegative constants $a$ and $\beta \leq 1$ such that
\begin{equation}\label{hdemaccr}
Re<Tx, Ax> + a \left\|x\right\|^2+\beta \left\|Tx\right\|^2 \geq 0 ,\quad  \text{ for all }  x\in \mathcal{D}  (T).
\end{equation}
If $\beta<1$, then  $T+A$ is m-accretive and also   the closure of $T+A$ is m-accretive for $\beta=1$, \cite[Theorem 4.1]{Okazawa73}. Note that this result cover the case  of relatively bounded perturbation, see \cite[Remark 4.4]{Okazawa73}. There are many papers on the question of such perturbation, see \cite{Okazawa77, Okazawa82, Okazawa2002, Sohr81, Yoshikawa72} for more results.
%%%%%%%%%%%%%%%%%%%%%%%%%%%%%%%%%%%%%%%%%%%%%%%%%%%%%%%%%%%%%

The aim of  this paper is to establish a new perturbation results on the m-accretivity of the operator $T+A$. This may be viewed as a slight improvement and  generalization of the perturbation results, particularly, those of  Okazawa, \cite{Okazawa77, Okazawa73}. The following lemma is our partial answer to the question above.
%%%%%%%%%%%%%%%%%%%%%%%%%%%%%%%%%%%%%%%%%%%%%%%%%%%%%%%%%%%
\begin{lemma}\label{perurbation} Let $T$ and $A$ two  operators such that   $\mathcal{D} (T)\subset \mathcal{D} (A)$. Assume that $T$ is m-accretive, $A$ is accretive and  there exists $c\geq 0$, such that 
\begin{equation}\label{bmaccr}
Re<Tx, Ax> \geq c\left\|Ax\right\|^2  ,\quad  \text{ for all }  x\in \mathcal{D}  (T).
\end{equation}
If we take $b=\max \{c\geq 0 : \eqref{bmaccr}  \text{ holds } \}$, we have,
\begin{enumerate}
	\item if $0 \leq b \leq 1$, then   $T+ A$ is also m-accretive,
	\item if $b>1$ then $T+ A$ is m-$\omega$-accretive,  with $\omega= \pi/2-\arcsin(\dfrac{b-1}{b})$.
\end{enumerate}
\end{lemma}
%%%%%%%%%%%%%%%%%%%%%%%%%%%%%%%%%%%%%%%%%%%%%%%%%
Here, $T$ is m-$\omega$-accretive if $e^{\pm i\theta} T$ is  m-accretive  for $\theta=\frac{\pi}{2}-\omega$, $0 < \omega\leq  \pi/2$. In this case, $-T$ generates  an holomorphic  contraction semigroup on the sector $\left| arg(\lambda)\right|<\omega$. In this
connection, we note that  for any  $\varepsilon>0$
$$ \left\| (\lambda + T)^{-1} \right\|\leq   \dfrac{M_\varepsilon}{\left| \lambda\right| }, \qquad \text{ for } \left| arg(\lambda)\right|\leq \dfrac{\pi}{2}+\omega-\varepsilon $$
with $M_\varepsilon$ is independent of $\lambda$ (see \cite[pp. 490]{Kato}).

The novelty of the lemma is the optimality of $b$ such that \eqref{bmaccr} holds. Clearly,  \eqref{bmaccr} implies  $Re<Tx, Ax>\geq 0$ for all   $x\in \mathcal{D}  (T)$, this exactly the assumption of \cite[Theorem 2.]{Okazawa69}.  Hence, we conclude that  $T+ A$ is also m-accretive. Our result is a refinement of this result by given a more precise sector containing the numerical range in function of the constant  $b$.  Also,  from \eqref{bmaccr}, we have for $b>0$, 
\begin{equation}\label{bmaccr1}
\left\|Ax\right\| \leq \dfrac{1}{b}\left\|Tx\right\| ,\quad  \text{ for all }  x\in \mathcal{D}  (T). 
\end{equation}
Thus the assumption  \eqref{bmaccr}  is stronger than the relative boundedness with respect to $T$. In particular, if $b>1$ the lower bound $\dfrac{1}{b}<1$, so according to  \cite[Theorem 2.]{Gustafson66},  $T+  A$ is m-accretive. Here, we say more,  $T+ A$ is m-$\omega$-accretive with $\omega$ depends of  the lower bound   $\dfrac{1}{b}<1$.       

%%%%%%%%%%%%%%%%%%%%%%%%%%%%%%%%%%%%%%%%%%%%%%%%%%%%%%
\section{Proof of the Lemma}
%%%%%%%%%%%%%%%%%%%%%%%%%%%%%%%%%%%%%%%%%%%%%%%%%%%%%%%%%
\begin{proof}[Proof of Lemma \ref{perurbation}]Let  $b=\max \{c\geq 0 : \eqref{bmaccr}  \text{ holds } \}$. 
If $b=0$, this exactly the \cite[Theorem 2.]{Okazawa69}. Assume that $0\leq b\leq 1$. We obtain from \eqref{bmaccr}
\begin{align*}
0&\leq Re<Tx, Ax>-b \left\|Ax \right\|^2\\
& \leq Re<Tx, Ax>+ (\alpha - b) \left\|Ax \right\|^2\\
\end{align*}
for some $ \alpha> 1$. Using  \eqref{hdemaccr}, we get
$$0\leq Re<Tx, Ax>+ \dfrac{\alpha- b}{b^2}  \left\|Tx \right\|^2.$$
Choosing $\alpha$  such that   $\beta=\dfrac{\alpha- b}{b^2} <1$, by \eqref{hdemaccr} we conclude that $T+  A$ is m-accretive (cf.\cite[Theorem 4.1]{Okazawa73}).

Now, suppose that that $b>1$. Let $x\in \mathcal{D}  (T)$, then for every $t> 0$, we have
\begin{align*}
Re<tx+ Tx, Ax>& =tRe<x, Ax>+Re<Tx, Ax>\\
& \geq b \left\|Ax \right\|^2.
\end{align*}
Thus we have
	\begin{equation}\label{bmaccr01}
\left\|Ax \right\| \leq \dfrac{1}{b}\left\|tx+ Tx\right\|.
	\end{equation}
Since $T$ is m-accretive, then
$$
\left\|A(t+ T)^{-1}x \right\| \leq \dfrac{1}{b}\left\|x\right\| ,\quad  \text{ for all }  x\in \mathcal{H}.
$$	
Hence it follows that
\begin{equation}\label{bmaccr02}
\left\| A(t+ T)^{-1} \right\| \leq \dfrac{1 }{b}<1.
\end{equation}	
 Then the operator $I+ A(t+ T)^{-1}$ is invertible and
$$
\left\|(I+A(t+ T)^{-1})^{-1} \right\| \leq \dfrac{b }{b-1}.
$$
The fact that
$$t+ T+ A =[I+ A(t+ T)^{-1}](t+ T), $$
it follows that $-t\in \rho(T+A)$ and
$$
\left\|t(t+ T+A)^{-1} \right\| \leq \dfrac{b }{b-1}=M ,\quad  \text{ for all }  t>0, 
$$
with $M>1$. Since $T+A$ is accretive, $\rho(T+A)$ contains  also the half plane $\{z\in \mathbb{C}  : Re(z)<0 \} $. Put $S =\{z\in \mathbb{C} : z\neq 0; \left|arg(z) \right| <\pi/2-\arcsin(\dfrac{1}{M})=\theta \}$ and  $M':= 1/ \sin(\pi/2 -\theta')$ with $0<\theta<\theta'<\pi/2 $, clearly $M'>M$. Let $\mu \in \mathbb{C} $ such that $ \left|arg(\mu)  \right| \leq  \theta'$ and fix $\lambda$ with $Re \lambda=-t<0$. Let $\left|  \mu -\lambda \right|\leq \dfrac{\left|\lambda \right|}{M'}$, we have that $\left\|  (\mu -\lambda)(t + T+ A)^{-1}\right\|\leq \dfrac{M}{M'} <1.$ Hence it follows that $\mu \in \rho(T+A)$ and
$$ (\mu + T+ A)^{-1}  = (\lambda + T+ A)^{-1} [I+ (\mu -\lambda)(\lambda + T+ A)^{-1}]^{-1}. $$
Thus
\begin{align*}
\left\|\mu  (\mu + T+ A)^{-1} \right\| & \leq \dfrac{\left|\mu \right|}{\left|\lambda \right|}\dfrac{1}{1-\dfrac{M}{M'} }M \\
& \leq (1+\dfrac{1}{M'})\dfrac{1}{1-\dfrac{M}{M'} }M.
\end{align*}
On the other hand,
\begin{align*}
(1+\dfrac{1}{M'})\dfrac{1}{1-\dfrac{M}{M'} }M &=\dfrac{1+\sin(\pi/2 -\theta')}{ \sin(\pi/2 -\omega)-\sin(\pi/2 -\theta')}\\
&\leq \dfrac{1}{ \sin((\theta' -\theta)/2)\sin((\theta' +\theta)/2)}\\
&\leq \dfrac{1}{ \sin(\theta' -\theta) \sin(\theta)}\\
&\leq \dfrac{1}{ \sin(\theta' -\theta) \sin(\pi/2-\arcsin(\dfrac{1}{M}))}\\
&\leq \dfrac{1}{ \sin(\theta' -\theta) \cos(\arcsin(\dfrac{1}{M}))}\\
&\leq \dfrac{1}{ \sin(\theta' -\theta) \sqrt{1-\dfrac{1}{M^2}}}\\
&\leq \dfrac{M}{ \sin(\theta' -\theta)\sqrt{M^2-1}}.
\end{align*}
This implies that 
$$ \left\| (\mu + T+ A)^{-1} \right\|\leq  \dfrac{M}{\left|\mu \right|\sin(\theta' -\theta)\sqrt{M^2-1}}.$$ 
This shows that the sector $S$  belongs to  $\rho(T+A)$ and for any $\varepsilon>0$,
$$ \left\| (\mu + T+ A)^{-1} \right\|\leq  \dfrac{M_\varepsilon}{\left|\mu \right|} \quad \text{ for } \quad \left|arg(\mu) \right| \leq\pi/2-\arcsin(\dfrac{1}{M})+\varepsilon,$$
with $M_\varepsilon=\dfrac{M}{\sin(\varepsilon)\sqrt{M^2-1}}$ and $\theta' -\theta= \varepsilon$.  Clearly, $M_\varepsilon$ is independent of $\mu$.  Hence, $T+ A$ is m-$\omega$-accretive,  with $ \omega=  \pi/2-\arcsin(\dfrac{b-1}{b})$.
\end{proof}
%%%%%%%%%%%%%%%%%%%%%%%%%%%%%%%%%%%%%%%%%%%%%%%%%%%%%%%
%%%%%%%%%%%%%%%%%%%%%%%%%%%%%%%%%%%%%%%%%%%%%%%%%%%%%%%%%%%
\begin{remark}\label{remlem}
	\begin{enumerate}	 
	\item As seen in the last paragraph of the  proof,  the condition \eqref{hdemaccr} implies \eqref{bmaccr} at least for  $0\leq b\leq 1$. Thus  \cite[Theorem 4.1]{Okazawa73} is covered by Lemma \ref{perurbation}. 
	\item If the assumptions of  Lemma \ref{perurbation} are satisfied, we can see that $Re<tx+ Tx, Ax>\geq 0$ for all $x\in  \mathcal{D}  (T)$. Therefore  $ A(t+ T)^{-1}$ is bounded accretive operator for any $t>0$.
	\end{enumerate}
\end{remark}
%%%%%%%%%%%%%%%%%%%%%%%%%%%%%%%%%%%%%%%%%%%%%%%%%%%%%%%
\begin{corollary}Let $T$ and $A$ as in Lemma \ref{perurbation}  obeying \eqref{bmaccr}. Then
\begin{enumerate}
	\item $-(T+ A)$ generates contractive one-parameter semigroup for  $0 \leq b \leq 1$.
	\item $-(T+ A)$ generates contractive holomorphic one-parameter semigroup  with angle $\omega= \arcsin(\dfrac{b-1}{b})$ for $b>1$.
\end{enumerate}	
\end{corollary}

%%%%%%%%%%%%%%%%%%%%%%%%%%%%%%%%%%%%%%%%%%%%%%%%%%%%%%%%
\section{An application}
%%%%%%%%%%%%%%%%%%%%%%%%%%%%%%%%%%%%%%%%%%%%%%%%%%%%%%%%

 One of interest is the operator-norm error bound estimate for the exponential Trotter-Kato
product formula in the case of accretive perturbations, see \cite{Cachia2002, Kato74, Kato78} and \cite{Neidhardt2018} for a short survey. Let $A$ be  a semibounded from below densely defined self-adjoint operator and  $B$ an m-accretive operator in a Hilbert space $ \mathcal{H}$. 

In \cite[Theorem 3.4]{Cachia2002}  it has been  shown that if $B$ is  $A$-bounded with lower bound $<1$ and 
\begin{equation}\label{cnd}
\mathcal{D} ((A+B)^{\alpha})\subset \mathcal{D} (A^{\alpha})\cap \mathcal{D} ((B^*)^{\alpha})\neq \{0\} \quad \text{ for some } \alpha \in \left( 0. 1\right],  
\end{equation}
then there is a constant $L_{\alpha}>0$ such that the estimates
\begin{equation}\label{EST1}
\left\|\left(e^{-tB/n} e^{-tA/n}\right)^n -e^{-t(A+B)}  \right\| \leq L_{\alpha}\dfrac{\ln n}{n^\alpha} 
\end{equation}
and
\begin{equation}\label{EST2}
\left\|\left( e^{-tA^*/n}e^{-tB^*/n}\right)^n -e^{-t(A+B)^*}  \right\| \leq L_{\alpha}\dfrac{\ln n}{n^\alpha} 
\end{equation} 
hold for some $\alpha \in \left( 0. 1\right]$ and $n=1,2,\ldots$ uniformly in $t\geq0$. Here $T^\alpha$ denotes the fractional powers of an m-accretive operator, see \cite{Kato60, Kato61}.  

The aim of the present result is to extend  \cite[Theorem 3.4]{Cachia2002}. This
extension is accomplished by replacing the relative boundedness by the assumption  \eqref{bmaccr}.  More precisely, we have

%%%%%%%%%%%%%%%%%%%%%%%%%%%%%%%%%%%%%%%%%%%%%
\begin{theorem}\label{thmTK} Let $A$ be  a semibounded from below densely defined self-adjoint operator and  $B$ an m-accretive operator with \eqref{bmaccr} for some $b>1$. Assume that    \eqref{cnd} holds. Then there is a constant $L_{\alpha}>0$ such that the estimates \eqref{EST1} and \eqref{EST2} hold for some $\alpha \in \left( 0. 1\right]$ and $n=1,2,\ldots$ uniformly in $t\geq0$.
\end{theorem}
%%%%%%%%%%%%%%%%%%%%%%%%%%%%%%%%%%%%%%%%%%%%%%%%
\begin{proof}From \eqref{bmaccr}, we have for $b>1$, 
	\begin{equation}\label{bmaccr12}
	\left\|Bx\right\| \leq a\left\|Ax\right\| ,\quad  \text{ for all }  x\in \mathcal{D}  (A), 
	\end{equation}
with $a=\frac{1}{b}<1$.  Hence $B$ is  $A$-bounded with lower bound $a<1$. Also, by lemma \ref{perurbation},  $A+B$ is m-$\omega$-accretive,  with $\omega= \pi/2-\arcsin(\dfrac{b-1}{b})$. Now, all assumptions of \cite[Theorem 3.4]{Cachia2002}
are fulfilled. Hence we obtain the desired result.	
\end{proof}	
%%%%%%%%%%%%%%%%%%%%%%%%%%%%%%%%%%%%%%%%%%	
\begin{remark} It well known that, for an m-accretive operator $T$,  the fractional powers  $T^{\alpha}$ are m-$(\alpha\pi)/2$-accretive  and, if $\alpha\!\in\! (0,1/2)$, then
	$\mathcal{D} (T^{\alpha})= \mathcal{D} (T^{*\alpha})$, see \cite[Theorem 1.1]{Kato61}.  Since $A$, $B$ and $A+B$ are m-accretive operators, we deduce that
	\begin{equation*}
\mathcal{D} ((A+B)^{*\alpha})=	\mathcal{D} ((A+B)^{\alpha})\subset \mathcal{D} (A^{\alpha})\cap \mathcal{D} (B^{\alpha})=\mathcal{D} (A^{\alpha})\cap \mathcal{D} ((B^*)^{\alpha}), 
	\end{equation*}
	for some $\alpha \in \left( 0,  1/2\right[$. Thus, the condition \eqref{cnd} may be omitted in Theorem \ref{thmTK}   if we take   $\alpha \in \left( 0,  1/2\right[$ (cf. \cite[Theorem 4.1]{Cachia2002}).
\end{remark}

%%%%%%%%%%%%%%%%%%%%%%%%%%%%%%%%%%%%%%%%%%%%%%%%%%%%%%%%.

%%%%%%%%%%%%%%%%%%%%%%%%%%%%%%%%%%%%%%%%%%%%%%%%%%%%%%%%

%%%%%%%%%%%%%%%%%%%%%%%%%%%%%%%%%%%%%%%%%%%%%%%%%
\end{document}